\documentclass[10pt,a4paper]{amsart}
\usepackage[utf8]{inputenc}
\usepackage[T1]{fontenc}
\usepackage[english]{babel}
\usepackage{amsmath,amsfonts,amssymb,amsthm}
\usepackage{stmaryrd} 
\usepackage{array}
\usepackage{enumerate}
\usepackage{hyperref}
\usepackage{mathabx}
 \usepackage{amsbsy}
\usepackage{lineno}
\newcommand{\Q}{\mathbb{Q}}
\newcommand{\N}{\mathbb{N}}
\newcommand{\Z}{\mathbb{Z}}
\newcommand{\R}{\mathbb{R}}
\newcommand{\Gal}{\mathrm{Gal}}
\newcommand{\tors}{\mathrm{tors}}
\newcommand{\poubelle}[1]{}

\theoremstyle{plain}
\newtheorem{thma}{Theorem}[section]
\newtheorem{propa}[thma]{Proposition}
\newtheorem{cora}[thma]{Corollary}
\newtheorem{lmma}[thma]{Lemma}

\theoremstyle{plain}
\newtheorem{thm}{Theorem}[subsection]
\newtheorem{prop}[thm]{Proposition}

\newtheorem{lmm}[thm]{Lemma}
\newtheorem{conj}[thm]{Conjecture}
\newtheorem{rqu}[thm]{Remark}

\begin{document}
%\linenumbers

\title[Small points on elliptic curves]{Location of small points on an elliptic curve by an equidistribution argument}
\author{Arnaud Plessis}
\address{Arnaud Plessis : Academy of Mathematics ans Systems Science, Morningside Center of Mathematics, Chinese Academy of Sciences, Beijing 100190}
\email{plessis@amss.ac.cn}

\begin{abstract}
Let $E$ be an elliptic curve defined over a number field $K$ without complex multiplication. 
If $\Gamma \subset E(\overline{K})$ is a subgroup of finite rank, a very special case of a conjecture of R\'emond predicts that points of small height in $E(K(\Gamma))$ lie in the division group of $\Gamma$. 
Using an equidistribution argument, we will show that this conjecture is true for groups of rank arbitrarily large.
\end{abstract}

\maketitle 

\section{Introduction}
\subsection{Main result}
Throughout this article $\overline{K}$ denotes an algebraic closure of a field $K$. 
Consider for this introduction an elliptic curve $E$ defined over a number field $K$.
Let $\hat{h} : E(\overline{K}) \to \R$ denote the N\'eron-Tate height on $E$. 

The rank of a subgroup $\Gamma\subset E(\overline{K})$ is given by the maximal number of linearly independent elements in $\Gamma$. 
The groups of rank $0$ coincide with the subgroups of $E_\tors$, the group of torsion points of $E$.  
Then we define its division group by 
\[
\Gamma_\mathrm{div} = \{P\in E(\overline{K}) \; \vert \;  \exists n\geq 1, [n]P \in \Gamma\},
\]
where $[n] : E(\overline{K}) \to E(\overline{K})$ is the multiplication-by-$n$ map. 
Note that $\Gamma$ and $\Gamma_\mathrm{div}$ have the same rank. 

Let $L\subset \overline{K}$ be a field, and let $F\subset E(\overline{K})$ be a set.
Write $L(F)$ for the smallest subfield of $\overline{K}$ containing $L$ over which all elements of $F$ are rational. 
Finally, throughout this text, the mention "small points of $E(L)$ lie in $F$" is shorthand for "there exists a positive constant $c$ such that $\hat{h}(P)\geq c$ for all $P\in E(L)\backslash F$".

We focus in this paper on a special case of a conjecture due to R\'emond, see \cite[Conjecture 3.4]{Remond} for the general statement. 
\begin{conj} \label{conj 1}
Assume that $E$ does not have complex multiplication. 
Let $\Gamma\subset E(\overline{K})$ be a subgroup of finite rank.
\begin{enumerate} [(i)]
\item (strong form): There is a positive constant $c$ such that 
\[
\hat{h}(P) \geq \frac{c}{[K(\Gamma)(P) : K(\Gamma)]} \; \text{for all} \; P\in E(\overline{K})\backslash \Gamma_\mathrm{div}.
\]
\item (weak form): For any $\varepsilon >0$, there exists a positive constant $c_\varepsilon$ such that
\[
\hat{h}(P) \geq \frac{c_\varepsilon}{[K(\Gamma)(P) : K(\Gamma)]^{1+\varepsilon}} \; \text{for all} \; P\in E(\overline{K})\backslash \Gamma_\mathrm{div}.
\]
\item (degree one form): Small points of $E(K(\Gamma))$ lie in $\Gamma_\mathrm{div}$.
\end{enumerate}
\end{conj}
We clearly have $(i) \Rightarrow (ii) \Rightarrow (iii)$.
Denote by $O$ the neutral element of $E$.
The strong form generalizes Lehmer's problem (which we recover when $\Gamma=\{O\}$) as well as the relative Lehmer's problem (corresponding to the case $\Gamma= E_\tors$).

The degree one form was proved to be true in the following situations: 
\begin{itemize}
\item (Northcott) $\Gamma$ is finitely generated \cite{Northcott};
\item (Habegger) $K=\Q$ and $\Gamma=E_\tors$ \cite[Theorem 2]{Habegger}.
\end{itemize} 
Up to my knowledge, these are the only two known cases going in the direction of Conjecture \ref{conj 1}.  

Let $S$ be a set of rational primes, and let $\Gamma\subset E(\overline{K})$ be a subgroup. 
We define the $S$-division group of $\Gamma$ as the subgroup $\Gamma_{\mathrm{div}, S}$ consisting of points $P\in E(\overline{K})$ for which there is $n\in\N=\{1,2,\dots\}$, whose prime factors are in $S$, such that $[n]P\in \Gamma$.

Several authors provided examples of fields $L\subset \overline{K}$ for which small points of $E(L)$ lie in $E_\tors$ \cite{Zhang, Baker, BakerSilverman, BakerPetsche, Pottmeyer, Sahu, Plessis} (in addition to \cite{Northcott, Habegger} cited above). 
However, up to my knowledge, no example has been provided of a field $L\subset \overline{K}$ for which it is possible to precisely locate small points of $E(L)$ if the latter do not lie in $E_\tors$.

Let $p$ be a rational prime, and let $\Gamma$ be the group generated by a single point $R\in E(\overline{K})\backslash E_\tors$. 
Clearly, $\hat{h}([1/p^n]R)=\hat{h}(R)/p^{2n}\to 0$ and $[1/p^n]R\in \Gamma_{\mathrm{div}, \{p\}}$ for all $n\in\N$. 
Hence, small points of $E(K(\Gamma_{\mathrm{div}, \{p\}}))$ do not lie in $E_\tors$. 
This simple example highlights the difficulty of Conjecture \ref{conj 1} and also explains why it is largely open. 
The main purpose of this paper is, under some conditions, to precisely locate small points of $E(K(\Gamma_{\mathrm{div}, \{p\}}))$.

If $v$ is a finite place of $K$, we denote by $K^{ur,v}$ the maximal algebraic extension of $K$ which is unramified above $v$. 
Next, if $N\in\N$ is an integer, we write $\zeta_N$ for a root of unity of order $N$. 
A more precise statement of the theorem below will be given in Theorem \ref{thm 4}. 

\begin{thm} \label{thm 2}
Let $E$ be an elliptic curve defined over a number field $K$ with split multiplicative reduction at a finite place $v$ of $K$.
Take a torsion-free subgroup $\Gamma\subset E(K)$ and assume that $\zeta_p\in K$, where $p$ denotes the rational prime below $v$.
Then there is an integer $c_0$ satisfying the following: For any set $S$ of rational primes greater than $c_0$, small points of $E(K^{ur,v}(\Gamma_{\mathrm{div}, S}))$ lie in $\Gamma_\mathrm{div}$. 
\end{thm}

\begin{rqu}
\rm{$a)$ This theorem is an elliptic analogue of \cite[Theorem 1.3]{Amoroso}, see also \cite[Th\'eor\`eme 1.8]{Plessis}. 

$b)$ The case where $S$ is the empty set is due to Baker \cite[Section 5, Case 1]{Baker}, see also \cite[Theorem 1.1]{Pottmeyer} for a statement explicitly written. 
Note that this case, and therefore our theorem, does not work anymore if $E$ has potential good reduction at $v$ \cite[Proposition 5.6]{Pottmeyer}.

$c)$ The assumption $\zeta_p\in K$ is not restrictive to study Conjecture \ref{conj 1} since we can enlarge $K$ as we wish.
In return, $S$ can contain more rational primes. 
In other words, without this assumption, we can prove the same result, but with a larger $c_0$.

$d)$ Throughout this text $\langle X\rangle$ denotes the group generated by a subset $X$ of a group $G$. 
The conditions relating to $\Gamma$ are not restrictive too. 
Let us see why. 
In Conjecture \ref{conj 1}, $\Gamma\subset E(\overline{K})$ is a subgroup of finite rank. 
It means that there exist linearly independent points $R_1,\dots, R_r\in E(\overline{K})\backslash E_\tors$ such that $\Gamma_{\mathrm{div}}$ is the set of $T+\sum_{j=1}^r [q_j]R_j$ with $T\in E_\tors$ and $q_j\in\Q$. 
It clearly suffices to prove Conjecture \ref{conj 1} with the larger subgroup $\Gamma_{\mathrm{div}}=\langle R_1,\dots, R_r\rangle_{\mathrm{div}}$ to deduce it in full generality. 
As we can also enlarge $K$, we can reduce the study of Conjecture \ref{conj 1} to the case that $R_1,\dots,R_r\in E(K)\backslash E_\tors$, that is $\langle R_1,\dots,R_r\rangle \subset E(K)$ is torsion-free. 

$e)$ Our proof does not provide us an effective lower bound for the N\'eron-Tate height of points in $E(K^{ur,v}(\Gamma_{\mathrm{div}, S}))\backslash E_\tors$.
The reason is that we will use an equidistribution theorem due to Chambert-Loir (see Theorem \ref{thm 3}) of which a quantitative version seems unknown. 
If it was known, then the arguments of this paper should permit us to make Theorem \ref{thm 2} effective.}
\end{rqu}

\subsection{Sketch of the proof of Theorem \ref{thm 2}} \label{subsection 1.2}
Throughout this text, a sequence of points $(P_n)_n$ in $E(\overline{K})$ is said to be a sequence of small points if $\hat{h}(P_n)\to 0$. 

Suppose that $(P_n)_n$ is a sequence of small points in $E(K^{ur,v}(\Gamma_{\mathrm{div}, S}))$.

$\mathbf{Step \; 1}$ (Proposition \ref{prop 1}): A result of Chambert-Loir on equidistribution of small points on elliptic curves with bad reduction allows us to construct, for all $n$ large enough, a subset $X_n\subset \Gal(\overline{K}/K)$ such that $\sigma P_n - P_n \in E_\tors$ for all $\sigma\in X_n$.
 
$\mathbf{Step \; 2}$ (after Lemma \ref{lmm 9}): Step 1 permits us to put in place a descent argument (Proposition \ref{prop 2}).
 This descent is then used to construct a finite extension $L/K^{ur,v}$ and a sequence of points $(T_n)_n$ in $\Gamma_\mathrm{div}$ such that $P_n-T_n\in E(L)$ for all $n$. 

$\mathbf{Step \; 3}$ (Lemma \ref{lmm 9}): We work with a new height $\hat{h}_\Gamma$, which is invariant under translation of points in $\Gamma_\mathrm{div}$ and satisfies $0\leq \hat{h}_\Gamma(P)\leq \hat{h}(P)$ for all $P\in E(\overline{K})$. 
Combining a theorem of Baker with a result due to R\'emond (Lemma \ref{lmm 7}), we show that the height $\hat{h}_\Gamma(P_n-T_n)$ cannot be arbitrarily small if $P_n-T_n\in E(L)\backslash \Gamma_\mathrm{div}$. 
But $0\leq \hat{h}_\Gamma(P_n-T_n)= \hat{h}_\Gamma(P_n) \leq \hat{h}(P_n) \to 0$, which leads to $P_n-T_n\in \Gamma_\mathrm{div}$, that is, $P_n \in \Gamma_\mathrm{div}$, for all $n$ large enough. 

\section{Preliminaries}
Let $K$ be a number field, and let $v$ be a finite place of $K$. 
We denote by $\vert . \vert_v$ the normalized $v$-adic absolute value on $K$, that is, $\vert p \vert_v = p^{-1}$, where $p$ is the rational prime below $v$.
The completion of $K$ with respect to $\vert. \vert_v$ is denoted with $K_v$. 
We use the same symbol $\vert . \vert_v$ for the absolute value on $K_v$, which we then extend to $\overline{K_v}$.
\subsection{Algebraic number theory}
Let $K$ be a number field, and let $v$ be a finite place of $K$.
We refer to \cite[Chapter II, § 2, Corollaire 3]{Serre} for the lemma below.
\begin{lmm} \label{lmm 1}
Two conjugate elements over $K_v$ have the same $v$- adic absolute value. 
\end{lmm}
We now fix a field embedding $\overline{K} \to \overline{K_v}$. 
We can thus see $\overline{K}$ as a subfield of $\overline{K_v}$.
Consider for this subsection a finite Galois extension $L/K$. 
For each place $w$ of $L$ above $v$, there exists a field embedding $\iota_w : L \to L_w$ such that $\vert \iota_w(.)\vert_v = \vert . \vert_w$. 
Moreover, $\iota_w$ is the identity on $K$ and $\iota_w(L)$ is dense in $L_w$ \cite[Proposition 1.3.1]{BombieriGubler}.
\begin{lmm} \label{lmm 2}
 The field embeddings $\tau\iota_w : L\to \overline{K_v}$ are pairwise distinct, where $w$ ranges over all places of $L$ above $v$ and $\tau$ runs over all elements of $\Gal(L_w/K_v)$. 
Moreover, each field embedding from $L$ to $\overline{K_v}$ fixing $K$ is of this form.  
\end{lmm}
\begin{proof}
Let $w_i$ be (with $i=1,2$) a place of $L$ above $v$, and let $\tau_i\in \Gal(L_{w_i}/K_v)$ such that $\tau_1 \iota_{w_1} = \tau_2 \iota_{w_2}$. 
Let $i\in\{1,2\}$. 
 Lemma \ref{lmm 1} implies $\vert \tau_i\iota_{w_i}(.)\vert_v =\vert \iota_{w_i}(.)\vert_v=\vert . \vert_{w_i}$.
Thus $\vert . \vert_{w_1} = \vert . \vert_{w_2}$.
Hence ${w_1}=w_2$, and so $\tau_1=\tau_2$ on $\iota_{w_1}(L)$. 

Again, Lemma \ref{lmm 1} gives $\vert \tau_i x - \tau_i y\vert_v = \vert x-y \vert_v$ for all $x,y\in L_{w_1}$.
Consequently, $\tau_i$ is uniformly continuous on $L_{w_1}$. 
The first part of the lemma arises from the continuous extension theorem since $\tau_1=\tau_2$ on $\iota_{w_1}(L)$ and $\iota_{w_1}(L)$ is dense in $L_{w_1}$.

As $K\subset \overline{K_v}$, the number of field embeddings $L\to \overline{K_v}$ which are the identity on $K$ is $[L : K]$. 
 We just constructed $\sum_w[L_w : K_v ]$ of them, where $w$ runs over all places of $L$ above $v$. 
The second part of this lemma follows since this sum is equal to $[L :K]$ by \cite[Corollary 1.3.2]{BombieriGubler}.
\end{proof} 
For any finite place $w$ of $L$, the decomposition group of the extension $w\vert w_{\vert K}$ is denoted by $D(w\vert w_{\vert K}) = \{\psi\in \Gal(L/K), \psi w = w\}$.

\begin{lmm} \label{lmm 3}
Let $w$ be a place of $L$ above $v$, and let $\sigma\in D(w \vert v)$. 
Then there exists $\sigma_w\in \Gal(L_w/K_v)$ such that $\iota_w \sigma = \sigma_w \iota_w$.
\end{lmm}
\begin{proof}
By Lemma \ref{lmm 3}, we have $\iota_w \sigma = \tau \iota_{w'}$ for some place $w'$ of $L$ above $v$ and some $\tau\in \Gal(L_{w'}/K_v)$. 
Then $\vert \sigma(.)\vert_w = \vert . \vert_{\sigma^{-1} w} = \vert . \vert_w$ since $\sigma\in D(w\vert v)$; whence $\vert \iota_w \sigma(.)\vert_v = \vert . \vert_w$.
Lemma gives us \ref{lmm 1} gives $\vert \tau\iota_{w'}(.)\vert_v = \vert \iota_{w'}(.)\vert_v = \vert . \vert_{w'}$.
So $\vert . \vert_w = \vert . \vert_{w'}$, which leads to $w=w'$.
The lemma follows since $\tau\in \Gal(L_{w'}/K_v)$.
\end{proof}

\subsection{Tate curves}
Fix once and for all an elliptic curve $E$ defined over a number field $K$ with split multiplicative reduction at a finite place $v$ of $K$.
We assume that $\zeta_p\in K$, where $p$ denotes the rational prime below $v$.
This hypothesis will only be used in the proof of Lemma \ref{lmm 5}. 
The objects $E, K, v$ and $p$ are as in Theorem \ref{thm 2}. 

The Tate uniformization provides a surjective morphism of $\Gal(\overline{K_v}/K_v)$-modules $\phi : \overline{K_v}^* \to E(\overline{K_v})$ with kernel $q^\Z =\{q^n,n\in\Z\}$ for some $q\in K_v^*$ of $v$-adic absolute value $<1$.
In particular, $\phi(L^*)= E(L)$ for each algebraic extension $L/K_v$ \cite[C.14]{Silverman}. 
We also fix a field embedding $\overline{K}\to \overline{K_v}$, which allows us to see from now $\overline{K}$ as a subfield of $\overline{K_v}$. 

Say $N\in\N$ is an integer. 
We define $E[N]$ as the set of torsion points killed by $N$.
Then, for $a\in\overline{K_v}$, denote by $a^{1/N}$ for a choice of solution of $(a^{1/N})^N=a$. 
It is defined up to $N$th-root of unity. 

\begin{lmm} \label{lmm 4}
For all integers $N\in\N$, we have $E[N]=\langle \phi(\zeta_N), \phi(q^{1/N})\rangle$. 
\end{lmm}
 
\begin{proof}
Let $T=\phi(u)\in E(\overline{K_v})$, and let $N\in\N$. 
Then $T\in E[N]$ if and only if $u^N\in q^\Z$ if and only if $u\in \langle \zeta_N, q^{1/N}\rangle$ if and only if $T=\phi(u)\in \langle \phi(\zeta_N), \phi(q^{1/N})\rangle$. 
\end{proof}

\begin{lmm} \label{lmm4.5}
Let $u\in\overline{K_v}^*$. 
Then $K_v(\phi(u))=K_v(u)$. 
\end{lmm}

\begin{proof}
As $u\in K_v(u)^*$, we have $\phi(u)\in \phi(K_v(u)^*)=E(K_v(u))$, which leads to the inclusion $K_v(\phi(u))\subset K_v(u)$.  
Similarly, $\phi(u)\in E(K_v(\phi(u)))=\phi(K_v(\phi(u))^*)$. 
Thus $\phi(u)=\phi(w)$ for some $w\in K_v(\phi(u))$. 
It means that $u=wq^n$ for some $n\in\Z$. 
As $q\in K_v$, we conclude $u\in K_v(\phi(u))$ and the lemma follows. 
\end{proof}

\subsection{Equidistribution} \label{subsection 2.2} 
Let $\hat{h}: E(\overline{K})\to \R$ denote the N\'eron-Tate height.
It is non-negative, invariant under Galois conjugation and vanishes precisely at $E_\tors$. 
It is also quadratic, that is,
\[
\forall n\in \Z, \forall P\in E(\overline{K}), \; \hat{h}([n]P)= n^2\hat{h}(P).
\]
This implies the assertion 
\[ 
\forall P\in E(\overline{K}), \forall T\in E_\tors, \; \hat{h}(P+T)= \hat{h}(P).
\]
Finally, it also satisfies the parallelogram law, that is, 
\[
\forall P, Q\in E(\overline{K}), \; \hat{h}(P+Q) + \hat{h}(P-Q)= 2(\hat{h}(P)+\hat{h}(Q)).
\]
We refer to \cite[Chapter VIII.9]{Silverman} for more information on $\hat{h}$.

Recall that $\phi : \overline{K_v}^*\to E(\overline{K_v})$ is a surjective homomorphism with kernel $q^\Z$.
If $x, y \in \overline{K_v}^*$ are preimages under $\phi$ of a point $P\in E(\overline{K_v})$, then $y=q^n x$ for some $n\in\Z$. 
Hence, the coset $l_v(P):= \log \vert x\vert_v/\log\vert q\vert_v +\Z \in \R/\Z$ is a well-defined element. 
The topological group $\R/\Z$ is homeomorphic to the unit circle and thus equipped with the unique Haar measure $\mu_{\R/\Z}$ of total mass $1$.  
The next theorem is due to Chambert-Loir \cite[Corollaire 5.5]{ChambertLoir}, see also \cite[Theorem 5]{Habegger} for a statement more in line with the spirit of this article. 
\begin{thm} \label{thm 3}
Let $(P_n)_n$ be a sequence of small points in $E(\overline{K})\backslash E_\tors$.
For all $n$, let $L_n/K(P_n)$ be a finite extension. 
Then, for each continuous function $f: \R/\Z \to \R$, we have 
\begin{equation} \label{eq. 1}
\frac{1}{[L_n : K]} \sum_\psi f(l_v(\psi(P_n))) \to \int f \mu_{\R/\Z},
\end{equation}
where $\psi$ runs over all field embeddings $L_n \to \overline{K_v}$ which are the identity on $K$. 
\end{thm}
\begin{rqu}
\rm{The left-hand side in \eqref{eq. 1} does not depend on $L_n$. 
This comes from the fact that any field embedding $K(P_n) \to \overline{K_v}$ which is the identity on $K$ extends in $[L_n : K(P_n)]$ different ways to a field embedding $L_n \to \overline{K_v}$.} 
\end{rqu}

Throughout this article $\#F$ denotes the cardinality of a finite set $F$. 
For any Galois extension $L/M$ and any subgroup $H$ of $\Gal(L/M)$, we write $L^H$ for the subfield of $L$ fixed by $H$. 
\begin{prop} \label{prop 1}
Let $(Q_n)_n$ be a sequence of small points in $E(\overline{K})$.
For all $n$, choose a finite Galois extension $L_n/K$ containing $K(Q_n)$. 
Let $H_n$ be a normal subgroup of $\Gal(L_n/K)$.
If the sequence $(\#H_n)_n$ is bounded, then we have, for all $n$ large enough, $\sigma Q_n - Q_n\in E_\tors$ for all $\sigma\in X_n= \bigcup_{w\vert v} D(w\vert w \cap L_n^{H_n})\subset H_n$, where $w$ ranges over all places of $L_n$ above $v$.
\end{prop}
\begin{proof}
Assume by contradiction that the conclusion of this proposition does not hold. 
It means that there exist infinitely many $n$ such that $\sigma_n Q_n - Q_n \notin E_\tors$ for some $\sigma_n\in X_n$ (that is, for all $n\in\N$ by taking a suitable subsequence).

Let $n\in\N$.
For brevity, put $L=L_n, \sigma=\sigma_n$ and $Q=Q_n$.
Set \[ C= \{\Psi\in \Gal(L/K), \Psi\sigma\Psi^{-1} = \sigma\}\] the centralizer of $\sigma$. 
By definition of $X_n$, there exists a place $\nu$ of $L$ above $v$ such that $\sigma\in D(\nu\vert \nu_{\vert L^{H_n}})$.
Write $C\nu=\{\Psi \nu, \Psi\in C\}$ and \[ \mathcal{F}=\{\tau\iota_w : L \to \overline{K_v}, w\in C \nu, \tau\in \Gal(L_w/K_v)\}.\]
As $L/K$ is Galois, Lemma \ref{lmm 2} gives $\#\mathcal{F}=\sum_{w\in C \nu} [L_w : K_v] = (\# C \nu)[L_{\nu} : K_v]$.
Moreover, the number of places of $L$ above $v$ is $[L : K]/[L_{\nu} : K_v]$.
The group $\Gal(L/K)$ acts transitively (for the natural action) on all places of $L$ above $v$. 
Hence, the cardinality of the orbit $C \nu$ is bounded from below by
\[
\#(C\nu) \geq \frac{\# C}{\#\Gal(L/K)} \frac{[L:K]}{[L_{\nu} : K_v]}.
\]
Let $N$ be an upper bound of the sequence $(\# H_n)_n$, which is bounded by assumption. 
As $H_n\ni \sigma$ is a normal subgroup of $\Gal(L/K)$, the orbit of $\sigma$ under the action by conjugation of the group $\Gal(L/K)$ is included in $H_n$. 
As $C$ is the stabilizer of $\sigma$ for this action, the orbit-stabilizer theorem proves that $\#\Gal(L/K)/\#C \leq \#H_n\leq N$.
Combining all of the above, we conclude that $\mathcal{F}$ has cardinality at least $[L:K]/N$.

Let $\psi=\tau\iota_w\in \mathcal{F}$ with $w\in C\nu$ and $\tau\in \Gal(L_w/K_v)$. 
Then $w = \Psi \nu$ for some $\Psi\in C$. 
Thus $\sigma = \Psi\sigma\Psi^{-1}$ and we easily check that $\Psi\sigma\Psi^{-1} w =w$; whence $\sigma\in D(w\vert v)$.
Lemma \ref{lmm 3} provides an element $\sigma_w\in \Gal(L_w/K_v)$ such that $\iota_w \sigma=\sigma_w \iota_w$. 
Recall that $\phi : \overline{K_v}^* \to E(\overline{K_v})$ is a surjective morphism of $\Gal(\overline{K_v}/K_v)$-modules. 
Let $u\in \overline{K_v}^*$ such that $\iota_w Q = \phi(u)$. 
A short calculation gives
\[
\psi (\sigma Q-Q) = \tau\iota_w\sigma Q - \tau\iota_w Q = \tau\sigma_w \phi(u)-\tau\phi(u) = \phi((\tau \sigma_w u)/\tau u).
\]
Lemma \ref{lmm 1} implies $\vert (\tau \sigma_w) u\vert_v = \vert \tau u\vert_v = \vert u \vert_v$, and so $l_v(\psi (\sigma Q - Q))=\Z$.
In conclusion, for each non-negative continuous function $f : \R/\Z\to \R^+$, we have
\begin{equation} \label{eq. 2} 
\sum_\psi \frac{f(l_v(\psi(\sigma Q-Q)))}{[L:K]} \geq \sum_{\psi\in \mathcal{F}} \frac{f(l_v(\psi(\sigma Q-Q)))}{[L:K]}= \frac{f(\Z)\#\mathcal{F}}{[L : K]} \geq \frac{f(\Z)}{N}, 
\end{equation}
where in the left-hand side, $\psi$ ranges over all field embeddings $L\to \overline{K_v}$ which are the identity on $K$. 
Some properties of $\hat{h}$ recalled earlier in this section lead to
\[
0\leq \hat{h}(\sigma Q - Q)\leq 2(\hat{h}(\sigma Q) + \hat{h}(Q)) = 4\hat{h}(Q)
\]
and the squeeze theorem proves that $\hat{h}(\sigma_n Q_n - Q_n) \to 0$. 
As $\sigma_n Q_n - Q_n\notin E_\tors$ for all $n$, we can then apply Theorem \ref{thm 3} to the sequence $(\sigma_n Q_n - Q_n)_n$. 
The left-hand side in \eqref{eq. 2} therefore goes to $\int f d\mu_{\R/\Z}$; whence $\int f d\mu_{\R/\Z}\geq f(\Z)/N$. 
The desired contradiction follows since there exist non-negative continuous functions $f : \R/\Z \to \R$ taking the value $1$ at $\Z$ and whose integral over $\R/\Z$ is $\leq 1/(2N)$.
\end{proof}

\section{A descent argument}
Recall that $\phi : \overline{K_v}^* \to E(\overline{K_v})$ is a surjective morphism of $\Gal(\overline{K_v}/K_v)$-modules with kernel $q^\Z$.
Furthermore, $\vert q\vert_v < 1$ and $\phi(L^*)=E(L)$ for each algebraic extension $L/K_v$. 
We also recall that $p$ is the rational prime below $v$ and that $\zeta_p\in K$. 
We only need this hypothesis to get the next lemma. 

Fix for this section an odd rational prime $l$ different to $p$. 
The goal of this section is to establish the descent argument that we will use to prove our theorem. 

\begin{lmma} \label{lmm 5}
Let $L/M$ be a Galois extension of local fields with $K_v\subset M$. 
Assume that $L/M$ is totally ramified of degree $l$. 
Let $\sigma\in \Gal(L/M)$, and let $Q\in E(L)$ such that $\sigma Q-Q\in E_\tors$. 
Then $\sigma Q-Q \in  \phi(\zeta_l)\Z$.
\end{lmma}
\begin{proof}
Let $u\in L^*$ such that $Q=\phi(u)$. 
Then
\begin{equation} \label{eq. 3}
\sigma Q - Q = \sigma \phi(u) - \phi(u) = \phi ((\sigma u)/u).
\end{equation}
Let $n\in\N$ be the order of $\sigma Q - Q$. 
By \eqref{eq. 3}, $n$ is the smallest positive integer for which $((\sigma u)/u)^n$ lies in the kernel of $\phi$, that is, $q^\Z$. 
On the one hand, $\vert q \vert_v <1$. 
On the other hand, $\vert (\sigma u)/u \vert_v = 1$ by Lemma \ref{lmm 1}. 
In conclusion, $((\sigma u)/u)^n =1$.
By the foregoing, $(\sigma u)/u$ has order $n$. 
Thanks to \eqref{eq. 3}, it suffices to get $n\vert l$ to deduce the lemma.
We have $\sigma u = \zeta_n u$ for some $\zeta_n\in L$. 

Put $n=m p^d$ with $p\nmid m$. 
Then $\zeta_n = \zeta_m \zeta_{p^d}$ for some $\zeta_m, \zeta_{p^d}\in L$.
 By \cite[Chapter IV, §4, Proposition 16]{Serre}, the extension $M(\zeta_m)/M$ is unramified since $p\nmid m$. 
The extension $L/M$ being totally ramified by hypothesis, we conclude $\zeta_m\in M$. 

Let $k\in\N$ such that $\sigma \zeta_{p^d}= \zeta_{p^d}^k$. 
As $\sigma u = \zeta_m \zeta_{p^d} u$ and $\sigma \zeta_m = \zeta_m$ (because $\zeta_m\in M$), we get, with the help of an easy induction,
\[
\sigma^i u = \zeta_m^i \zeta_{p^d}^{k^{i-1}+k^{i-2}+\dots+1} u
\]
for all $i\in \N$. 
Taking $i=l$, we obtain $\zeta_m^l \zeta_{p^d}^{k^{l-1}+\dots+1} =1$ since $\Gal(L/M)$ has order $l$ by hypothesis. 
It means, since $p$ and $m$ are coprime, that
\begin{equation} \label{eq. 4}
m\vert l \; \; \text{and} \; \; k^{l-1}+\dots+1 \equiv 0 \; (p^d).
\end{equation}
If $d=0$, then $n=m$ and the lemma follows from \eqref{eq. 4}. 
We now assume that $d\geq 1$ and we want to arrive at a contradiction. 
As $\zeta_{p^d}\in L$, the multiplicativity formula for degrees implies that $[M(\zeta_{p^d}) : M]$ divides $l=[L:M]$. 
But it also divides $[\Q_p(\zeta_{p^d}) : \Q_p(\zeta_p)]=p^{d-1}$ since $\zeta_p\in K\subset M$. 
So $[M(\zeta_{p^d}) : M]=1$ because $l\neq p$ by assumption. 
In particular, $\sigma \zeta_{p^d}=\zeta_{p^d}$; whence $k\equiv 1 \; (p^d)$. 
Combining this congruence and \eqref{eq. 4}, we get $l\equiv 0 \; (p^d)$, which is absurd since $d\geq 1$ and $l\neq p$. 
\end{proof}

Until the end of this text, denote by $O$ the neutral element of $E$.

The cyclic group $\phi(\zeta_l)\Z$ has cardinality either $1$ or $l$ since $[l]\phi(\zeta_l)=\phi(1)=O$. 
It cannot be trivial since otherwise, we would get $\zeta_l\in q^\Z$. 
As $\vert \zeta_l\vert_v=1$ and $\vert q\vert_v<1$, it follows that $\zeta_l=1$, which is absurd.
Hence $\phi(\zeta_l)\Z$ has cardinality $l$.

Given a finite extension $L/K$, a set $Y=\{P_1,\dots, P_s\}\subset E(\overline{K})$ and a point $Q\in E(L(Y))$, our descent argument consists in finding (under some conditions on $L, Y$ and $Q$) integers $n_1,\dots, n_s\in\{0,\dots, l-1\}$ such that $Q-\sum_{i=1}^s [n_i]T_i \in E(L)$. 

The lemma below corresponds to the base case $s=1$, the general situation (see Proposition \ref{prop 2}) arising from an easy induction on $s$. 
\begin{lmma} \label{lmm 6}
Let $L/K(E[l])$ be a finite extension, let $P\in E(\overline{K})$ and let $Q\in E(L(P))$. 
Assume that: 
\begin{enumerate} [(i)]
\item $L(P)/L$ is Galois and its Galois group is isomorphic to $(\Z/l\Z)^m$ for some $m\in\{0,1,2\}$; 
\item $\sigma Q-Q, \sigma P-P\in E_\tors $ for all $\sigma\in \bigcup_{w\vert v} D(w\vert w\cap L)$, where $w$ ranges over all places of $L(P)$ above $v$;
\item There is a place $w$ of $L(P)$ above $v$ such that $L(P)_w/L_w$ is totally ramified of degree $l$; 
\item If $m=2$, then there exists a subextension $M/K$ of $L/K$ satisfying: 
\begin{enumerate} [(a)]
\item $L/M$ and $L(P)/M$ are Galois; 
\item There exist $\psi\in\Gal(L(P)/M)$ and $\sigma\in D(w\vert w\cap L)$ such that the groups $\langle \psi^k\sigma\psi^{-k}\rangle$ are pairwise distinct when $k\in\{1,\dots,l\}$. 
\end{enumerate}
\end{enumerate}
Then $Q-[n]P\in E(L)$ for some integer $n\in\{0,\dots,l-1\}$. 
\end{lmma}

\begin{proof}
\underline{Case $m=0$:} Obvious since then $Q\in E(L(P))=E(L)$.

\underline{Case $m=1$:} Gathering together $(i)$ and $(iii)$ proves that $\Gal(L(P)/L)$ is generated by an element $\sigma$ belonging to $D(w\vert w\cap L)$.
Lemma \ref{lmm 3} claims that $\sigma_w\iota_w=\iota_w\sigma$ for some $\sigma_w\in \Gal(L(P)_w/L_w)$. 
Consequently, \[\iota_w(\sigma Q-Q)=\sigma_w(\iota_w Q)-(\iota_w Q)\in E_\tors \; \text{and} \; \iota_w(\sigma P-P)=\sigma_w(\iota_w P)-(\iota_w P)\in E_\tors  \]
(these points are torsion by $(ii)$). 
Thanks to $(iii)$, it arises from Lemma \ref{lmm 5} applied to $M=L_w, L=L(P)_w, \sigma = \sigma_w$ and $Q\in\{\iota_w Q, \iota_w P\}$ that $\sigma Q-Q$ and $\sigma P-P$ belong to $\iota_w^{-1} \phi(\zeta_l)\Z$. 
This group being cyclic of order $l$, it follows that $\sigma Q-Q = [n](\sigma P-P)$ for some $n\in\{0,\dots,l-1\}$ because $\sigma P\neq P$. 
The point $Q-[n]P$ is therefore fixed by $\sigma$, that is, $Q-[n]P\in E(L(P)^{\langle \sigma \rangle})=E(L)$.

\underline{Case $m=2$:} Put $\sigma_k=\psi^k\sigma\psi^{-k}$ and $w_k=\psi^k w$. 
The fact that $L/M$ is Galois by $(iva)$ shows that $\sigma_k\in \psi^k D(w\vert w\cap L)\psi^{-k}=D(w_k \vert w_k \cap L)$. 
As both $L(P)/M$ and $L/M$ are Galois by $(iv a)$, it follows that $L(P)_{w_k}/L_{w_k}=L(P)_w/L_w$. 
This extension is totally ramified of degree $l$ by $(iii)$. 
Following the lines of the previous case, we get $\sigma_k Q-Q, \sigma_k P-P\in T_k\Z$ for all $k\in\{1,\dots,l\}$, where $T_k=\iota_{w_k}^{-1}\phi(\zeta_l)$.

\underline{Subcase 1: $\sigma_k Q\neq Q$ for all $k\in\{1,\dots,l\}$.}  Let $k\in\{1,\dots,l\}$. 
As $\sigma_k P \neq P$, we deduce that $\sigma_k Q-Q, \sigma_k P-P \in (T_k\Z)\backslash \{O\}$. 
The cyclicity of $T_k\Z$ allows us to construct $n_k\in\{1,\dots,l-1\}$ such that $\sigma Q-Q= [n_k](\sigma P-P)$. 
This leads to $Q-[n_k]P\in E(L(P)^{\langle \sigma_k\rangle})$. 
But $k\in\{1,\dots,l\}$ and $n_k\in\{1,\dots,l-1\}$. 
The pigeonhole principle guarantees us the existence of two distinct integers $k,k'\in\{1,\dots,l\}$ satisfying $n:=n_k=n_{k'}$. 
Gathering together $(i)$ and $(iv b)$, we conclude that $\Gal(L(P)/L)$ is generated by $\sigma_k$ and $\sigma_{k'}$; whence \[Q-[n]P\in E\left(L(P)^{\langle \sigma_k\rangle} \cap L(P)^{\langle \sigma_{k'}\rangle}\right)=E\left(L(P)^{\langle \sigma_k, \sigma_{k'}\rangle}\right)=E(L). \]

\underline{Subcase 2: There is $k\in\{1,\dots,l\}$ such that $\sigma_k Q=Q$.} For brevity, assume that $k=1$. 
The application $f : \Gal(L(P)/L) \to E(\overline{K}), \tau  \mapsto  \tau P-P$ is injective since the identity is the unique element of $\Gal(L(P)/L)$ fixing $P$. 
By the foregoing, we have $\sigma_2 Q-Q, f(\sigma_2)\in T_2\Z$. 
As $T_2\in E[l]\subset L$ is fixed by $\sigma_2$, we get, thanks to a small calculation, $\sigma_2^nQ-Q, f(\sigma_2^n)\in T_2\Z$ for all $n\in\{1,\dots,l-1\}$. 
By $(i)$ and $(iv b)$, the group $\Gal(L(P)/L)$ is generated by $\sigma_1$ and $\sigma_2$. 
As $l$ is odd, it follows that $\sigma_3=\sigma_2^{l_2}\sigma_1^{l_1}$ for some $l_1,l_2\in\{1,\dots,l-1\}$ (we have $l_1,l_2\neq 0$ by $(iv b)$). 
But $\sigma_3 Q-Q=\sigma_2^{l_2}Q - Q$ since $\sigma_1$ fixes $Q$; whence $\sigma_3 Q- Q\in T_2\Z \cap T_3\Z$. 

If $T_3\Z=T_2\Z$, then $f(\sigma_3)\in T_2\Z$ by the foregoing. 
As $\sigma_3, \sigma_2,\dots, \sigma_2^{l-1}$ are $l$ pairwise distinct non-trivial elements of $\Gal(L(P)/L)$, the injectivity of $f$ shows that the $l$ points $f(\sigma_3), f(\sigma_2), \dots, f(\sigma_2^{l-1})$ are pairwise distinct and different to $O$. 
This is impossible since they all lie in $(T_2\Z)\backslash \{O\}$, which has cardinality $l-1$. 
So $T_3\Z\neq T_2\Z$. 
But then, $\sigma_3 Q - Q= O$ since $T_2\Z$ and $T_3\Z$ are two distinct groups with order $l$. 
Combining $(i)$ and $(iv b)$, we obtain that $\Gal(L(P)/L)$ is generated by $\sigma_1$ and $\sigma_3$. 
These two homomorphisms fixing $Q$, we conclude $Q\in E\left(L(P)^{\langle \sigma_1, \sigma_3\rangle}\right) = E(L)$, which ends the proof of the lemma.   
\end{proof}

We can now state our descent argument whose proof consists in iterating several times the previous lemma. 

\begin{propa} \label{prop 2}
Let $L/K(E[l])$ be a finite extension, let $\{P_1,\dots, P_s\}\subset E(\overline{K})$, and let $Q\in E(L(P_1,\dots,P_s))$. 
Assume that for all $j\in\{1,\dots,s\}$, we have: 
\begin{enumerate} [(i)]
\item $L(P_1,\dots,P_j)/L$ is Galois and the Galois group of $L(P_1,\dots,P_j)/L(P_1,\dots,P_{j-1})$ is isomorphic to $(\Z/l\Z)^{m_j}$ for some $m_j\in\{0,1,2\}$; 
\item $\sigma Q-Q, \sigma P_j-P_j\in E_\tors $ for all $\sigma\in \bigcup_{w\vert v} D(w\vert w\cap L)$, where $w$ ranges over all places of $L(P_1,\dots,P_s)$ above $v$; 
\item There is a place $w_j$ of $L(P_1,\dots,P_j)$ above $v$ such that $L(P_1,\dots,P_j)_{w_j}/L(P_1,\dots,P_{j-1})_{w_j}$ is totally ramified of degree $l$; 
\item If $m_j=2$, then there is a subextension $M_j/K$ of $L(P_1,\dots,P_{j-1})/K$ such that: 
\begin{enumerate} [(a)]
\item $L(P_1,\dots,P_{j-1})/M_j$ and $L(P_1,\dots,P_j)/M_j$ are Galois; 
\item There exist $\psi_j\in\Gal(L(P_1,\dots,P_j)/M_j)$ and $\sigma_j\in D(w_j\vert w_j\cap L(P_1,\dots,P_{j-1}))$ such that the groups $\langle \psi_j^k\sigma_j\psi_j^{-k}\rangle$ are pairwise distinct when $k\in\{1,\dots,l\}$. 
\end{enumerate}
\end{enumerate}
Then $Q-\sum_{j=1}^s [n_j]P_j\in E(L)$ for some integers $n_1,\dots,n_s\in\{0,\dots,l-1\}$. 
\end{propa}

\begin{proof}
By induction on $s$. 
The base case $s=1$ corresponds to Lemma \ref{lmm 6}. 
It is obvious that we can apply Lemma \ref{lmm 6} to $L=L(P_1,\dots, P_{s-1})$ and $P=P_s$ (take $j=s$ in $(i)-(iv)$ above). 
So we can find an integer $n_s\in\{0,\dots,l-1\}$ such that $Q-[n_s]P_s\in E(L(P_1,\dots,P_{s-1}))$. 
We now want apply our induction to the set $\{P_1,\dots,P_{s-1}\}$ and to $Q=Q-[n_s]P_s$. 
It is clear that $(i), (iii)$ and $(iv)$ hold since the larger set $\{P_1,\dots,P_s\}$ satisfies them by assumption. 
Prove $(ii)$.
Let $w$ be a place of $L(P_1,\dots,P_{s-1})$ above $v$, and let $\sigma\in D(w\vert w\cap L)$. 
As $L(P_1,\dots, P_s)/L$ and $L(P_1, \dots, P_{s-1})/L$ are Galois by $(i)$, \cite[Chapter 1, § 7, Proposition 22 $b)$]{Serre} tells us that $\sigma$ can be extended to an element $\tilde{\sigma}\in D(w'\vert w\cap L)$, where $w'$ is any place of $L(P_1,\dots,P_s)$ above $w$. 
By $(ii)$, we have $\sigma P_j - P_j = \tilde{\sigma} P_j - P_j\in E_\tors$ for all $j\leq s-1$. 
Finally, a small calculation gives \[\sigma (Q-[n_s]P_s) - (Q-[n_s]P_s) = \tilde{\sigma} Q - Q -[n_s](\tilde{\sigma} P_s - P_s),\] which is therefore a torsion point according to $(ii)$. 
We can now applied our induction to the situation described above in order to get $Q-[n_s]P_s-\sum_{j=1}^{s-1}[n_j]P_j \in E(L)$ for some integers $n_1,\dots, n_{s-1}\in \{0,\dots,l-1\}$. 
This finishes the proof.
\end{proof}

\section{Proof of Theorem \ref{thm 2}}
Let $\Gamma\subset E(K)$ be a subgroup. 
For $P\in E(\overline{K})$, set 
\[
\hat{h}_\Gamma(P)=\inf\{\hat{h}(P+Q), Q\in \Gamma_\mathrm{div}\}.
\]
Clearly, $0\leq \hat{h}_\Gamma(P)\leq \hat{h}(P)$ and $\hat{h}_\Gamma$ is invariant under translation of points in $\Gamma_\mathrm{div}$. 
We recover the N\'eron-Tate height if $\Gamma\subset E_\tors$ since $\hat{h}$ is invariant under translation of torsion points.  

The lemma below is along the lines of \cite[Lemme 3.5]{Remond}.
It shows that the height $\hat{h}_\Gamma$ is particularly well adapted to study the degree one form of Conjecture \ref{conj 1}.  
Recall that $\Gamma$ and $\Gamma_\mathrm{div}$ have the same rank (see the introduction).

\begin{lmma} [R\'emond] \label{lmm 7}
Let $L/K$ be an algebraic extension, and let $\Gamma\subset E(K)$ be a subgroup.
The following assertions are equivalent: 
\begin{enumerate}[(i)]
\item small points of $E(L)$ lie in $\Gamma_{\mathrm{div}}$;
\item there is $c>0$ such that $\hat{h}_\Gamma(P)\geq c$ for all $P\in E(L)\backslash \Gamma_{\mathrm{div}}$. 
\end{enumerate}
\end{lmma}

\begin{proof}
Clearly, $(ii)\Rightarrow (i)$ since $\hat{h}_\Gamma(P)\leq \hat{h}(P)$.
As $\Gamma\subset E(K)$, we have the inclusions $\Gamma\subset \Gamma_{\mathrm{div}} \cap E(L)\subset \Gamma_{\mathrm{div}}$. 
Consequently, $\Gamma, \Gamma_{\mathrm{div}} \cap E(L)$ and $\Gamma_{\mathrm{div}}$ have the same rank, which is finite since $\Gamma\subset E(K)$. 
The implication $(i)\Rightarrow (ii)$ now arises from \cite[Corollary 2.3]{Pottmeyer2}.
\end{proof}

Recall that $\phi : \overline{K_v}^* \to E(\overline{K_v})$ is a surjective homomorphism with kernel $q^\Z$ and that we fixed a field embedding $\overline{K}\to \overline{K_v}$.
Thus each point of $E(\overline{K})$ can be expressed as $\phi(a)$ for some $a\in\overline{K_v}^*$, what we do from now. 
We now fix a torsion-free group $\Gamma\subset E(K)$ and a generating set $\{R_1,\dots,R_r\}$ of $\Gamma$.
Put $R_i=\phi(a_i)$.

Let $v(.)$ be the $v$-adic valuation on $K_v^*$.
Write $\mathcal{T}$ for the set of rational primes $l$ such that $\Gal(K(E[l^d])/K)\simeq \mathrm{Aut} E[l^d]$ for all $d\in\N$ and not dividing \[2 v(q) p \#((\Gamma_\mathrm{div} \cap E(K))/\Gamma).\] 
By Serre's open image theorem, $\mathcal{T}$ contains all but finitely many rational primes. 
We can now state the precise version of Theorem \ref{thm 2}. 

\begin{thma} \label{thm 4}
Let $S\subset \mathcal{T}$ be a finite subset.
Then small points of $E(K^{ur,v}(\Gamma_{\mathrm{div}, S}))$ lie in $\Gamma_\mathrm{div}$. 
\end{thma}

It is clear that $([v(q)]\Gamma)_\mathrm{div}=\Gamma_\mathrm{div}$ and that the index of $[v(q)]\Gamma$ in $\Gamma$ is $v(q)^r$ (because the quotient group is in bijection with the set of $\sum_{i=1}^r [n_i]R_i$, where $n_i\in\{0,\dots,v(q)-1\}$). 
Thus, the set $\mathcal{T}$ remains unchanged if we replace $\Gamma$ with $[v(q)]\Gamma$. 
Next, as each element of $S$ does not divide $v(q)$, B\'ezout's identity claims that $\Gamma_{\mathrm{div}, S}=([v(q)]\Gamma)_{\mathrm{div}, S}$. 
Hence, it suffices to prove Theorem \ref{thm 4} for the subgroup $[v(q)]\Gamma$. 
The interest of this reduction is that $[v(q)]R_i=\phi(a_i^{v(q)})=\phi(a_i^{v(q)}q^{-v(a_i)})$ and it is easy to check that the $v$-adic valuation of $a_i^{v(q)}q^{-v(a_i)}$ is $0$. 
In conclusion, we can assume without loss of generality that $v(a_i)=0$ for all $i\in\{1,\dots,r\}$.

Say $N\in\N$ is an integer.
Fix from now a choice of $\zeta_N,q^{1/N}, a_1^{1/N},\dots,a_r^{1/N}$ and put $[1/N]R_i=\phi(a_i^{1/N})$. 
By abuse of notation, we define $[M/N]R_i$ as the point $[M]([1/N]R_i)$. 
The field $K(E[N], [1/N]R_i)$ is Galois over $K$. 
As everything is now fixed, we ease notations and for any field $F\subset \overline{K}$, we put \[F(N)= F(E[N], [1/N]R_1,\dots, [1/N]R_r).\]
We now fix a finite set $S\subset \mathcal{T}$.
Let $P\in E(K^{ur,v}(\Gamma_{\mathrm{div},S}))$, and write $\Lambda_P$ for the set of integers $N\in\N$ satisfying:
\begin{enumerate} [(i)]
\item the prime factors of $N$ lie in $S$; 
\item there exists $T\in \Gamma_{\mathrm{div}}$ such that $P-T \in E(K^{ur,v}(N))$.
\end{enumerate}
Put $k=\prod_{l\in S} l$.
The group $\Gamma_{\mathrm{div}, S}$ is then the set of $Q\in E(\overline{K})$ such that $[k^n]Q\in\Gamma$ for some $n\geq 0$. 
Concretely, \[\Gamma_{\mathrm{div}, S}=\bigcup_{n\geq 0} \left\{T+\sum_{j=1}^m [a_j/k^n]R_j, T\in E[k^n], a_j\in\Z\right\}.\]
Hence $K(\Gamma_{\mathrm{div}, S})= \bigcup_{n\geq 0}K(k^n)$. 
Thus $P\in E(K^{ur,v}(k^n))$ for some $n\geq 0$, and so $k^n\in \Lambda_P$. 
Denote by $N_P$ the minimal element of $\Lambda_P$.
Finally, let $T_P\in \Gamma_{\mathrm{div}}$ denote a point such that $P-T_P\in E(K^{ur,v}(N_P))$.

We can now start the proof of Theorem \ref{thm 4}.
Let $(P_n)_n$ be a sequence of small points in $E(K^{ur,v}(\Gamma_{\mathrm{div}, S}))$. 
We want to show that $P_n\in \Gamma_{\mathrm{div}}$ for all $n$ large enough. 
For brevity, put $N_n=N_{P_n}$ and $T_n=T_{P_n}$. 
Let us start by treating the case where the sequence $(N_n)_n$ is bounded. 

\begin{lmma} \label{lmm 9}
If the sequence $(N_n)_n$ is bounded, then $P_n\in \Gamma_{\mathrm{div}}$ for all $n$ large enough.
\end{lmma}

\begin{proof}
Let $N$ be the least common multiple of $N_1, N_2,\dots$.
A result due to Baker \cite[Section 5, Case 1]{Baker} (see also \cite[Theorem 1.1]{Pottmeyer} for a reference explicitly written) asserts that small points of $E(K^{ur,v}(N))$ lie in  $E_\tors \subset \Gamma_{\mathrm{div}}$.
Lemma \ref{lmm 7} provides a positive constant $c$ such that $\hat{h}_\Gamma(Q)\geq c$ for all $Q\in E(K^{ur,v}(N))\backslash \Gamma_{\mathrm{div}}$. 
As $\hat{h}_\Gamma$ is invariant under translation of points in $\Gamma_\mathrm{div}$, we conclude \[ 0\leq \hat{h}_\Gamma(P_n-T_n)= \hat{h}_\Gamma(P_n)\leq \hat{h}(P_n) \] and the squeeze theorem proves that $\hat{h}_\Gamma(P_n-T_n) \to 0$. 
But $P_n-T_n\in E(K^{ur,v}(N))$ by definition of $T_n$ and $N$. 
In conclusion, $P_n-T_n\in \Gamma_{\mathrm{div}}$ for all $n$ large enough and the lemma follows since $T_n\in\Gamma_{\mathrm{div}}$. 
\end{proof}

Assume that the sequence $(N_n)_n$ is unbounded.
By Lemma \ref{lmm 9}, it suffices to get a contradiction to deduce Theorem \ref{thm 4}. 
By considering a subsequence of $(P_n)_n$ if needed, we can assume that $N_n \to +\infty$.
Our descent argument will show that $N_n$ cannot be the minimal element of $\Lambda_{P_n}$, thus deriving to the desired contradiction.
But first of all, note that $\hat{h}(P_n-T_n)$ can be arbitrarily large although $\hat{h}(P_n)\to 0$.
See below how to circumvent this problem. 

Let $n\in\N$. 
As $T_n\in \Gamma_{\mathrm{div}}$, there exist $T'_n\in E_\tors, x_{j,n}\in \Z$ and $t_n\in \N$ such that $T_n=T'_n+\sum_{j=1}^r [x_{j,n}/t_n]R_j$.
The field $\Q$ being archimedean, we get $x_{j,n}= q_{j,n} t_n N_n^{-1}+r_{j,n}$ for some $q_{j,n}\in \Z$ and some $r_{j,n}\in\Q$ satisfying $0\leq r_{j,n} < t_n N_n^{-1}$.
A short calculation gives 
\begin{equation} \label{eq def Q}
Q_n :=  P_n-T_n+ \sum_{j=1}^r [q_{j,n}/N_n] R_j =  P_n-T''_n-\sum_{j=1}^r [r_{j,n}/t_n]R_j
\end{equation}
 for some $T''_n\in E_\tors$. 
Note that $Q_n\in E(K^{ur,v}(N_n))$.
Using the fact that $\hat{h}$ is invariant under translation of torsion points, then the parallelogram law, we get \[0  \leq \hat{h}(Q_n) = \hat{h}\left(P_n-\sum_{j=1}^r [r_{j,n}/t_n]R_j\right)  \leq 2\hat{h}(P_n)+2\sum_{j=1}^r \hat{h}([r_{j,n}/t_n]R_j).\]
As $\hat{h}$ is quadratic and $r_{j,n}/t_n < N_n^{-1}$, we finally conclude \[0\leq \hat{h}(Q_n) \leq 2\hat{h}(P_n)+ 2 N_n^{-2} \sum_{j=1}^r \hat{h}(R_j),\]
and therefore $\hat{h}(Q_n)\to 0$ since $N_n\to +\infty$ and $\hat{h}(P_n)\to 0$.  

The fact that $S$ is finite and that $N_n\to +\infty$ permits us to find a rational prime $l\in S$ such that $d_n$, the $l$-adic valuation of $N_n$, goes to $+\infty$. 
Let $u$ be the largest integer such that $\zeta_{l^u}\in K_v(\zeta_l)$. 
Capelli's Lemma \cite[Chapter VI, Theorem 9.1]{Lang} claims that the extension $K_v(\zeta_{l^t})/K_v(\zeta_l)$ has degree $l^{t-u}$ for all $t\geq u$. 
By removing the first terms of the sequence $(N_n)_n$ if needed, we can assume that $d_n\geq \max\{2, u+1\}$ for all $n$. 
In particular, $l$ divides $N_n/l$. 

It is well known that $K^{ur,v}/K$ is Galois. 
So there is a number field $M_n\subset K^{ur,v}$, Galois over $K$, such that $Q_n\in E(M_n(N_n))$ for all $n$. 
In addition, as $K^{ur,v}$ contains all roots of unity with order coprime to $p$, we can enlarge $M_n$ so that it contains all $(p^{f(N_n^2!)}-1)$th-roots of unity for all $n$, where $f$ denotes the inertia degree of the extension $K_v/\Q_p$.
This reduction guarantees us that the unramified extension of degree $N_n^2!$ over $K_v$ is contained in $(M_n)_w$ for all places $w$ of $M_n$ above $v$ \cite[Chapter IV, §4, Corollaire 1]{Serre}.

We plan to apply Proposition \ref{prop 2} to $Q=Q_n$, $L=M_n(N_n/l), s=r+2$ and to \[P_s= \phi(q^{1/l^{d_n}}), P_{s-1}= \phi(\zeta_{l^{d_n}}) \; \; \text{and} \; \;  P_j= [1/l^{d_n}]R_j\] for all $j\in\{1,\dots,r\}$.
By Lemma \ref{lmm 4}, $(P_s, P_{s-1})$ forms a basis of $E[l^{d_n}]$. 
We can hence identify $P_s$ and $P_{s-1}$ with points in $E(\overline{K})$. 

Note that $Q\in E(M_n(N_n))= E(L(P_1,\dots,P_s))$ and that $E[l]\subset L$ since $l$ divides $N_n/l$.
Thus, if all five conditions of Proposition \ref{prop 2} are fulfilled for at least one $n$, then our descent would prove that \[Q_n-\sum_{j=1}^s [n_j]P_j \in E(M_n(N_n/l))\subset E(K^{ur,v}(N_n/l))\] for some integers $n_1,\dots,n_s$. 
Each $P_j$ belonging to $\Gamma_\mathrm{div}$, it follows from \eqref{eq def Q} that $P_n-T'''_n\in  E(K^{ur,v}(N_n/l))$ for some $T'''_n\in \Gamma_{\mathrm{div}}$. 
The prime factors of $N_n/l$ being in $S$, we get $N_n/l\in \Lambda_{P_n}$, thus contradicting the minimality of $N_n$. 

\underline{Proof of $(i)$:} Let $i\in\{1,\dots,s\}$. 
By definition of $d_n$, we have $l^{d_n-1}$ divides $N_n/l$. 
Hence $[l]P_i\in L$, and so the Galois conjugates of $P_i$ over $L$ belong to $P_i+E[l]$, which is included in $L(P_i)$ since $E[l]\subset E(L)$. 
Thus $L(P_i)/L$ is Galois and its Galois group is isomorphic to a subgroup of $E[l]\simeq (\Z/l\Z)^2$. 
This shows $(i)$. 

\underline{Proof of $(ii)$:} By the foregoing, we have $\sigma P_j - P_j \in E[l]$ for all $j\in\{1,\dots,s\}$ and all $\sigma\in \Gal(\overline{K}/L)$.  
Furthermore, $[L(P_1,\dots,P_s) : L]\leq l^{2s}$. 
As $\hat{h}(Q_n)\to 0$, Proposition \ref{prop 1} applied to $L_n=L(P_1,\dots,P_s)$ and $H_n=\Gal(L(P_1,\dots,P_s)/L)$ allows us to obtain $(ii)$ if $n$ is large enough. 
Choose from now such an integer $n$ and put $M=M_n, N=N_n$ and $d=d_n$. 

\underline{Proof of $(iii)$:}
We proceed by proving a series of results.

\begin{lmma} \label{lmm 5.5}
Let $K\subset F\subset F'$ be a tower of number fields. 
Take a finite place $w$ of $F'$ as well as $\tau\in\Gal(\overline{K}/K)$. 
Then the extensions $F'_w/F_w$ and $(\tau F')_{\tau w}/(\tau F)_{\tau w}$ have the same ramification index and the same inertia degree. 
\end{lmma}

\begin{proof}
Denote by $\mathcal{O}_k$ the ring of integers of a number field $k$. 
The isomorphism \[ \left(\mathcal{O}_{\tau F'}/\tau w\right)/\left(\mathcal{O}_{\tau F}/\tau w \cap \tau F\right) = \left((\tau \mathcal{O}_{F'})/\tau w\right)/ \left((\tau \mathcal{O}_F)/\tau (w\cap F)\right)  \simeq \left(\mathcal{O}_{F'}/w\right)/ \left(\mathcal{O}_F/(w\cap F)\right)  \] proves that $(\tau F')_{\tau w}/(\tau F)_{\tau w}$ and $F'_w/F_w$ have the same inertia degree. 
They also have the same ramification index since if $(w\cap F)\mathcal{O}_{F'}=w^e I$, where $e\in\N$ and $I$ is an ideal of $\mathcal{O}_{F'}$ coprime to $w$, then $(\tau w \cap \tau F)\mathcal{O}_{\tau F'}= (\tau w)^e \tau I$ with $\tau w $ and $\tau I$ coprime. 
This completes the lemma.
\end{proof}

Recall that $v(a_i)=0$ and that the prime divisors of $N$ all belong to $S$, which does not contain the prime factors of $p v(q)$ (in particular, $N$ and $v(q)$ are coprime).

\begin{lmma} \label{lmm102}
Let $w$ be a place of $M$ above $v$. 
Then $\zeta_N, a_1^{1/N},\dots,a_r^{1/N}\in M_w$.
\end{lmma}

\begin{proof}
Class field theory tells us that $K_v(\zeta_N)/K_v$ is unramified since $N$ is coprime to $p$. 
Next, Kummer theory claims that $K_v(\zeta_N, a_j^{1/N})/K_v(\zeta_N)$ is unramified because $v(a_j)=0$ and $N$ is coprime to $p$. 
In conclusion, $K_v(\zeta_N, a_j^{1/N})$ is an unramified extension of $K_v$ whose degree divides $N^2!$.  
The lemma now arises from the fact that $M_w$ contains the unramified extension of degree $N^2!$ over $K_v$.
\end{proof}

\begin{lmma} \label{lmm103}
Let $t\in\N$ be a divisor of $N$, and let $F/K_v(\zeta_t)$ be an unramified extension. 
Then $F(q^{1/t})/F$ is totally ramified of degree $t$. 
\end{lmma}

\begin{proof}
Let $w$ be the place of $K_v(\zeta_t)$, which is an unramified extension of $K_v$. 
By Kummer theory, $K_v(\zeta_t, q^{1/t})/K_v(\zeta_t)$ is totally ramified of degree $t$ since $N$, and so $t$, is coprime to $w(q)=v(q)$.
The lemma follows since $F/K_v(\zeta_t)$ is unramified.
\end{proof}

\begin{thma} \label{thm 5}
We have $\Gal(K(t)/K(E[t]))\simeq E[t]^r$ for all $t\in\N$ dividing $N$. 
\end{thma}

\begin{proof}
Applying \cite[Chapter V, Theorem 5.2]{Lang2} to $n=t$ and to $\{P_1,\dots,P_r\}=\{R_1,\dots,R_r\}$. 
All requested conditions are fulfilled by construction. 
\end{proof}
 
\begin{cora} \label{cor 1}
Let $t\in\N$ be an integer dividing $N$. 
Then $K([1/t]R_1,\dots,[1/t]R_r)$ and $K(E[t])$ are linearly disjoint over $K$. 
\end{cora}

\begin{proof}
Let $i\in\{1,\dots,r\}$. 
As $R_i\in E(K)$, the Galois conjugates of $[1/t]R_i$ over $K$ belong to $[1/t]R_i+E[t]$; whence $[K([1/t]R_i) : K]\leq t^2$. 
From Theorem \ref{thm 5}, we get \[t^{2r}=[K(t) : K(E[t])] \leq [K([1/t]R_1,\dots,[1/t]R_r) : K] \leq t^{2r}. \]
The corollary follows. 
\end{proof}

\begin{lmma} \label{lmm 105}
For all $j\in\{1,\dots,s\}$, there exists $\tau_j\in\Gal(\overline{K}/K)$ such that $\tau_j P_i=P_i$ for all $i<j$ and $\tau_j P_j = P_j+[\min\{1,s-j\}]P_s$. 
\end{lmma}

\begin{proof}
It is clear if $j=s$ because the identity works.
If $j\leq r=s-2$, then Theorem \ref{thm 5} applied to $t=l^d$ ensures us the existence of an element $\tau\in \Gal(\overline{K}/K)$ such that $\tau P_i = P_i + O $ if $i<j$ and $\tau P_i = P_i + P_s$ if $i=j$.
Finally, if $j=s-1$, then there is $\tau\in \Gal(K(E[l^d])/K)=\mathrm{Aut} E[l^d]$ such that $\tau P_{s-1}=P_{s-1}+P_s$. 
Corollary \ref{cor 1} applied to $t=l^d$ asserts that we can lift $\tau$ to an element $\tilde{\tau}\in \Gal(K(l^{d})/K(P_1,\dots,P_r))$ satisfying $ \tilde{\tau} P_{s-1}= \tau P_{s-1}$.
The lemma is now proved. 
\end{proof}

Let $j\in\{1,\dots,s\}$.
In order to unify the case $j\leq s-2$ and the case $j=s-1$, we put $\zeta_t=a_{s-1}^{1/t}$ for all $t\in\N$. 
From Lemma \ref{lmm 105}, we easily get $\tau_j P_j = \phi(x_j)$, where \[x_j= \begin{cases} 
 q^{1/l^d} \; \text{if} \; j=s \\ 
a_j^{1/l^d} q^{1/l^d} \;  \text{otherwise}
\end{cases}.\] 
Write $w_j$ for the place of $\tau_j(L(P_1,\dots,P_j))$ associated to the fixed embedding $\overline{K}\to \overline{K_v}$. 
Recall that $L=M(E[N/l], [l/N]R_1,\dots, [l/N]R_r)$ is Galois over $K$.
Thus 
\begin{multline*}
\tau_j(L(P_1,\dots,P_j))_{w_j} = L(P_1,\dots,P_{j-1}, \tau_j P_j)_{w_j} = \\  M_{w_j}(\phi(q^{l/N}), \phi(\zeta_{N/l}), \phi(a_1^{l/N}),\dots,\phi(a_r^{l/N}), \phi(a_1^{1/l^d}), \dots, \phi(a_{j-1}^{1/l^d}), \phi(x_j)).  
\end{multline*}
Combining Lemma \ref{lmm4.5} with Lemma \ref{lmm102} gives \[\tau_j(L(P_1,\dots,P_j))_{w_j} = M_{w_j}(q^{l/N}, q^{1/l^d})= M_{w_j}(q^{1/N}). \] 
Similarly, $\tau_j(L(P_1,\dots,P_{j-1}))_{w_j}=M_{w_j}(q^{l/N})$.  
Lemma \ref{lmm103} applied to $t=N$ and $F=M_{w_j}$ shows that the extension $\tau_j(L(P_1,\dots, P_j))_{w_j}/ \tau_j(L(P_1,\dots, P_{j-1}))_{w_j}$ is totally ramified of degree $l$. 
We can express $w_j$ as $w_j=\tau_j w'_j$ for some place $w'_j$ of $L(P_1,\dots,P_j)$ above $v$. 
Item $(iii)$ follows from Lemma \ref{lmm 5.5} applied to $F'=L(P_1,\dots, P_j), F= L(P_1,\dots, P_{j-1}), \tau=\tau_j$ and to $w=w'_j$. 

\underline{Proof of $(iv)$:} 
\underline{Case $j\geq r$:} Recall that $d\geq \max\{2, u+1\}$ and that $K_v(\zeta_{l^t}) / K_v(\zeta_l)$ has degree $l^{t-u}$ for all $t\geq u$. 
In particular, $\zeta_{l^d}\notin K_v(\zeta_{l^{d-1}})$.

\begin{lmma}
Let $T\in E[l^d]$, and let $F/M(E[l^{d-1}])$ be any extension.  
Then the degree of the extension $F(T) / F$ is either $1$ or $l$. 
\end{lmma}

\begin{proof}
As $[l]T\in E(F)$, the Galois conjugates of $T$ above $F$ belong to $T+E[l]$, which is included in $F(T)$ since $E[l]\subset E[l^{d-1}]\subset E(F)$. 
In particular, $F(T) / F$ is Galois and its Galois group is isomorphic to a subgroup of $(\Z/l\Z)^2$. 
Assume by contradiction that it is $(\Z/l\Z)^2$. 
Then \[l^2=[F(T) : F]\leq [K(E[l^{d-1}], T) : K(E[l^{d-1}])] \leq l^2\] and the restriction map $\Gal(F(T)/F) \to \Gal(K(E[l^{d-1}], T) / K(E[l^{d-1}]))$ is therefore an isomorphism. 
But $\zeta_{l^d}\in M\subset F$. 
So $\zeta_{l^d}\in K(E[l^{d-1}])$. 
The completion of $K(E[l^{d-1}])$ with respect to any place above $v$ is $K_v(\zeta_{l^{d-1}}, q^{1/l^{d-1}})$. 
By Lemma \ref{lmm103} applied to $t=l^{d-1}$ and to $F=K_v(\zeta_{l^{d-1}})$, the extension $K_v(\zeta_{l^{d-1}}, q^{1/l^{d-1}})/ K_v(\zeta_{l^{d-1}})$ is totally ramified. 
The extension $K_v(\zeta_{l^d})/K_v$ being unramified, we deduce that $\zeta_{l^d}\in K_v(\zeta_{l^{d-1}})$, which is absurd.
The lemma follows.
\end{proof}

Applying this lemma to $F=L(P_1,\dots, P_r)$ and to $T=P_{s-1}$, resp. $F=L(P_1,\dots,P_{s-1})$ and $P=T_s$, gives $m_{s-1}\leq 1$, resp. $m_s\leq 1$.

\underline{Case $j<r$:} 
Put $M_j=K([l]P_1,\dots,[l]P_j)$. 
We clearly have $[l]P_j\in E(L)$ and $E[l]\subset E(L)$. 
So the extensions $L(P_1,\dots,P_{j-1})/M_j$ and $L(P_1,\dots,P_j)/M_j$ are Galois, which proves $(iva)$. 
It remains to get $(ivb)$. 
For brevity, put $w=w'_j$ and take a generator $\sigma\in D(w\vert w\cap L(P_1,\dots,P_{j-1}))$.

Corollary \ref{cor 1} applied to $t=l^{d-1}$ shows that $K(E[l])$ and $M_j$ are linearly disjoint over $K$. 
Thus $\Gal(M_j(E[l])/M_j)\simeq \Gal(K(E[l])/K) = \mathrm{Aut} E[l]$. 
Let $(T,T')$ be a basis of $E[l]$.
As $[l]P_j\in E(L)$, there are two integers $a$ and $b$ satisfying \[\sigma P_j=P_j + [a]T+[b]T'. \]
As $\sigma$ does not fix $P_j$, we have either $a\not\equiv 0 \; (l)$ or $b\not\equiv 0 \; (l)$.
By considering the basis $(T',T)$ if needed, we can assume that $a\not\equiv 0 \; (l)$. 
Let $\psi\in\Gal(M_j(E[l])/M_j)$ denote the element satisfying $\psi T = T+T'$ and $\psi T' = T'$.
We have $\psi^k T = T +[k]T'$ and $\psi^k T' = T'$ for all $k\in\N$. 
We can extend $\psi$ to an element of $\Gal(L(P_1,\dots,P_j)/M_j)$; call it $\psi$ again.
As $[l]P_j\in M_j$, we can find two integers $c$ and $d$ satisfying \[\psi P_j = P_j + [c]T+[d]T'. \]
An easy induction gives us $\psi^kP_j = P_j+[kc]T+[kd+k(k-1)c/2]T'$ for all $k\in\N$. 

Assume by contradiction that the groups $\langle \psi \sigma \psi^{-1}\rangle, \dots, \langle \psi^l \sigma \psi^{-l}\rangle$ are not pairwise distinct. 
A small calculation shows that there is $k\in\{1,\dots,l-1\}$ such that $\psi^k\langle \sigma\rangle \psi^{-k} = \langle \sigma \rangle$.
Equivalently, there is $n\in\{1,\dots,l-1\}$ such that $\psi^k\sigma^n = \sigma\psi^k$. 
Thus, as $\sigma$ fixes $L$, and therefore $E[l]$, we get \[\sigma\psi^k P_j = P_j+[a]T+[b]T'+[kc]T+[kd+k(k-1)c/2]T'=\psi^k P_j + [a]T+[b]T'.\] 
A similar calculation proves that \[\psi^k\sigma^n P_j = \psi^k(P_j+[an]T+[bn]T')=\psi^k P_j+ [an]T+[ank]T'+[bn]T'.\] 
By identification, we get $an\equiv a \; (l)$ and $ank+bn\equiv b \; (l)$. 
As $a\not\equiv 0 \; (l)$, we get $n=1$, leading to $ak+b\equiv b \; (l)$, which is absurd since $a,k\not\equiv 0 \; (l)$. 
\qed

\subsection*{Acknowledgment}
I thank P. Habegger, L. Pottmeyer and G. R\'emond for answering my questions as well as the referees for his/her valuable suggestions. 
This work was funded by Morningside Center of Mathematics, CAS. 
\bibliographystyle{plain}

\end{document}